\documentclass[12pt]{amsart}

\usepackage{mathrsfs}

\usepackage{amsfonts, amssymb,amsmath, amsthm, amsxtra, latexsym, amscd}
\usepackage[latin1]{inputenc} 

\usepackage[english]{babel}

\usepackage{amsmath}

\usepackage{amsthm,amssymb,latexsym}

\usepackage{t1enc}

\usepackage{epic}

\usepackage{eepic}

\usepackage{xypic}

\usepackage{amsfonts, amssymb,amsmath, amsthm, amsxtra, latexsym, amscd}

\usepackage{graphics}

\title[Quaternionic discrete series for $Sp(1,1)$]{Quaternionic discrete series for $Sp(1,1)$}
\author{Henrik Seppänen}
\address{Department of Mathematics, Chalmers University of Technology 
and G\"oteborg University, G\"oteborg, Sweden}
\email{henriks@math.chalmers.se}

\newtheorem{prop}{Proposition}

\newtheorem{lemma}[prop]{Lemma}

\newtheorem{thm}[prop]{Theorem}

\theoremstyle{definition}

\theoremstyle{remark} \newtheorem*{remark}{Remark}

\begin{document}

\nocite{*}
\thispagestyle{empty}


\cleardoublepage

\pagestyle{headings}
\setcounter{page}{1}

\maketitle

\begin{abstract}
In this paper we study the analytic realisation of the discrete series representations for the group $G=Sp(1,1)$ as a subspace of the space of
square integrable sections in a homogeneous vector bundle over the symmetric space $G/K:=Sp(1,1) /(Sp(1) \times Sp(1))$.
We use the Szegö map to give expressions for the restrictions of the $K$-types occurring in the representation spaces to the submanifold 
$AK/K$. 
\end{abstract}

\section{Introduction}
In \cite{GW}, Gross and Wallach considered representations of simple Lie groups $G$ with maximal compact subgroup $K$ such that 
the associated symmetric space $G/K$ has a $G$-equivariant quaternionic structure (cf. \cite{wolf}).
This amounts to the group $K$ containing a normal subgroup isomorphic to $SU(2)$. In fact, there is an isomorphism
$K \cong SU(2) \times M$ for a subgroup $M \subseteq K$, and by setting $L:=U(1) \times M$, the associated homogeneous space
$G/L$ is fibred over $G/K$ with fibres diffeomorphic to $P^1(\mathbb{C})$.
The quaternionic discrete series representations are then realised on the Dolbeault cohomology groups $H^1(G/L, \mathcal{L})$, 
where $\mathcal{L} \rightarrow G/L$ is a 
holomorphic line bundle. In this model they are able to classify all the $K$-types occurring in each
of the obtained discrete series representations. Moreover, they consider the continuation of the discrete series and 
characterise the unitarisability of the underlying $(\mathfrak{g}, K)$-modules.

In this paper we consider another model of the quaternionic discrete series. If $\pi$ is a quaternionic discrete series
representation realised on the cohomology group $H^1(G/L, \mathcal{L})$, and $\tau$ is its minimal $K$-type, then the
\emph{Schmid D-operator} acts on the sections of the homogeneous vector bundle $G \times _K V_{\tau} \rightarrow G/K$ where
$V_{\tau}$ is some vector space on which the $K$-type is unitarily realised. The Hilbert space $\mbox{ker} D \cap L^2(G,\tau)$
then furnishes another realisation of the representation $\pi$. We consider the special case when $G=Sp(1,1)$. In this case
the symmetric space $G/K$ can be embedded into the bounded symmetric domain $SU(2,2)/S(U(2) \times U(2))$ consisting
of complex $2 \times 2$-matrices of norm less than one. The restriction of the Harish-Chandra embedding to $G/K$ then 
yields a global trivialisation of the vector bundle $G \times_K V_{\tau}$. In this model we compute
the restrictions to the submanifold $A \cdot 0$ \footnote{$A$ is associated with a particular Iwasawa decomposition
$G=NAK$.} of the highest weight vectors for the occurring $K$-types.
These functions turn out to be fibrewise highest weight vectors with a hypergeometric function as a coefficient.
Similar functions have been studied by Castro and Gr\"unbaum in \cite{castro-grunbaum}. Hypergeometric functions occur frequently
in representation theory, not only for Lie groups. For example, in \cite{opdam-hecke}, they play a role in the context of Hecke algebras.

We compute the $K$-types by using the Szeg\"o map defined by Knapp and Wallach in \cite{K-W} which exhibits any discrete series representation
as a quotient of a nonunitary principal series representations. The $K$-types are determined on the level of the principal series
representation, and then the Szeg\"o map is applied to compute the above mentioned restrictions.

The paper is organised as follows. In section 2 we explicitly state some results from the structure theory of the Lie group $Sp(1,1)$ 
that will be needed.
In section 3 we describe the models for the discrete series in the general context of induced representations, and also give an explicit global 
trivialisation. Section 4 describes the Szeg\"o map by Knapp and Wallach, and we also compute $K$-types on the level of 
a nonunitary principal series representation. In section 5 we compute the images of the $K$-types under the Szeg\"o map and 
trivialise them to yield vector valued functions. The main theorem of this paper is Theorem 8 of this section.

\section{Preliminaries}

\subsection{The quaternion algebra}
The quaternion algebra, $\mathbb{H}$, is a four-dimensional associative algebra over $\mathbb{R}$ with generators
$i,j,k$ satisfying the relations 
\begin{eqnarray*}
i^2=j^2=k^2=-1\\
ij=k, jk=i, ki=j \,\, \mbox{and}\\
ji=-ij, ik=-ki, kj=-jk.
\end{eqnarray*}
Moreover, $\mathbb{H}$ is equipped with an involution, $*$, given by
$$(a+bi+cj+dk)^*=a-bi-cj-dk, a,b,c,d \in \mathbb{R}.$$
The Euclidean norm on the vector space $\mathbb{R}^4 \simeq \mathbb{H}$ can be expressed in terms of this involution by
$$|(a,b,c,d)|^2=a^2+b^2+c^2+d^2=(a+bi+cj+dk)^*(a+bi+cj+dk).$$
It follows immediately that the quaternions of norm one, $Sp(1)$, form a group.
The algebra $\mathbb{H}$ can be embedded as a subalgebra of the algebra, $M_2(\mathbb{C})$, of $2 \times 2$ complex matrices by
\begin{eqnarray}
\iota : \mathbb{H} \rightarrow M_2(\mathbb{C}), \label{matrisid}
\end{eqnarray}
where 
\begin{eqnarray}
\iota(a+bi+cj+dk)=\left(\begin{array}{cc}
a+bi & c+di \\
-(c-di) & a-bi \\
\end{array}\right).
\end{eqnarray}
In particular, the generators $1,i,j,k$ are embedded as

\begin{eqnarray*}
\iota(1)=\left(\begin{array}{cc}
1 & 0 \\
0 & 1 \\
\end{array}\right), 
\iota(i)=\left(\begin{array}{cc}
i & 0 \\
0 & -i \\
\end{array}\right), 
\iota(j)=\left(\begin{array}{cc}
0 & 1 \\
-1 & 0 \\
\end{array}\right), 
\iota(k)=\left(\begin{array}{cc}
0 & i \\
i & 0 \\
\end{array}\right).
\end{eqnarray*}
The embedding $\iota$ also satisfies the relation 
$$\iota((a+bi+cj+dk)^*)=\left(\begin{array}{cc}
a-bi & -c-di \\
-(-c+di) & a+bi \\
\end{array}\right)\\
=\left(\begin{array}{cc}
a+bi & c+di \\
-(c-di) & a-bi \\
\end{array}\right)^*,$$
so $\iota$ is a homomorphism of involutive algebras.
We observe that, letting $z=a+bi, w=c+di$,  $$\iota(\mathbb{H})=\left\{\left(\begin{array}{cc}
z & w \\
-\overline{w} & \overline{z} \\
\end{array}\right)|\, z,w \in \mathbb{C}\right\}$$ and moreover, we have the identity
$$a^2+b^2+c^2+d^2=|z|^2+|w|^2= \det\left(\begin{array}{cc}
z & w \\
-\overline{w} & \overline{z} \\
\end{array}\right).$$ 
In particular, $$Sp(1) \simeq 
\left\{\left(\begin{array}{cc}
z & w \\
-\overline{w} & \overline{z} \\
\end{array}\right)|\,|z|^2+|w|^2=1\right\}=SU(2).$$

\subsection{The group Sp(1,1)}
The real vector space $\mathbb{H}^2 \cong \mathbb{R}^4$ is also equipped with the structure of an $\mathbb{H}$-module by
\begin{eqnarray}
\left(\alpha,(h_1,h_2)\right) \mapsto (\alpha h_1, \alpha h_2), \alpha, h_1, h_2 \in \mathbb{H}.
\end{eqnarray}
If we identify $\mathbb{H}^2$ with the set of $2 \times 1$ matrices over $\mathbb{H}$, 
there is a natural $\mathbb{H}$-linear action of the matrix group $GL(2,\mathbb{H})$ on $\mathbb{H}^2$ given by
\begin{eqnarray}
\left( \begin{array}{cc}  h_1 &  h_2 \end{array}\right) \mapsto \left( \begin{array} {cc} h_1 & h_2 \end{array}\right)
\left( \begin{array}{cc}
a & b\\
c & d\\
\end{array}\right).
\end{eqnarray}
Consider the real vector space $\mathbb{H}^2$ equipped with the nondegenerate indefinite Hermitian form 
\begin{eqnarray}
\langle \cdot, \cdot \rangle_{1,1}: ((h_1,h_2),(h_1',h_2')) \mapsto h_1(h_1')^*-h_2(h_2')^*.
\end{eqnarray}
Recall that the group $Sp(1,1)$ is defined as
\begin{eqnarray}
Sp(1,1):=\{g \in GL(2, \mathbb{H})|\langle gh, gh'\rangle_{1,1}=\langle h, h'\rangle_{1,1}\},
\end{eqnarray}
where $h:=(h_1,h_2), h':=(h_1',h_2')$.
The condition that the form $\langle, \rangle_{1,1}$ be preserved can be reformulated as
\begin{eqnarray}
g^*Jg=J, \label{matrisvillkor}
\end{eqnarray}
where $g^*=\left( \begin{array}{cc}
a & b\\
c & d \\
\end{array} \right)^*:=\left( \begin{array}{cc}
a^* & c^*\\
b^* & d^* \\
\end{array} \right)$,
and $J=\left( \begin{array}{cc}
1 & 0\\
0 & -1 \\
\end{array} \right)$.

The embedding \eqref{matrisid} induces an embedding (which we also denote by the same symbol)
\begin{eqnarray}
\iota: M_2(\mathbb{H}) \rightarrow M_4(\mathbb{C}) \label{matrisjota}
\end{eqnarray}
 by
\begin{eqnarray}
\left( \begin{array}{cc}
a & b\\
c & d \\
\end{array} \right) \mapsto \left( \begin{array}{cc}
\iota(a) & \iota(b)\\
\iota(c) & \iota(d) \\
\end{array} \right).
\end{eqnarray}
This embedding is a homomorphism of algebras with involution. Applying it to the identity \eqref{matrisvillkor} reveals that the image
of $Sp(1,1)$ is a subgroup of the group
\begin{eqnarray}
SU(2,2)&=&\{g \in M_4(\mathbb{C})|g^*\tilde{J}g=\tilde{J}, \mbox{det}g=1\},\\
\tilde{J}&:=&\left( \begin{array}{cc}
I_2 & 0\\
0 & -I_2 \\
\end{array} \right). \nonumber
\end{eqnarray}

\subsection{The symmetric space $B_1(\mathbb{H})=Sp(1,1)/(Sp(1) \times Sp(1))$}
Let $B_1(\mathbb{H})$ denote the unit ball
\begin{eqnarray}
B_1(\mathbb{H}):=\{h \in \mathbb{H}|\,|h|<1\}
\end{eqnarray}
in $\mathbb{H}$. The group $G:=Sp(1,1)$ acts transitively on $B_1(\mathbb{H})$ by the fractional linear action
\begin{eqnarray}
g(h):=(ah+b)(ch+d)^{-1}, \label{fraktionell}\\  
g=\left(\begin{array}{cc}
a & b\\
c & d\\
\end{array}\right) \in G, h \in B_1(\mathbb{H}). \nonumber
\end{eqnarray}
The isotropic subgroup for the origin is
\begin{eqnarray}
K:=G_0=\left\{ \left( \begin{array}{cc}
a & 0\\
0 & d\\
\end{array}\right) \in G\right\} \cong Sp(1) \times Sp(1),
\end{eqnarray}
and hence we have the description
\begin{eqnarray}
B_1(\mathbb{H})\cong G/K
\end{eqnarray}
of $B_1(\mathbb{H})$ as a homogeneous space.
Moreover, from eq. \eqref{matrisvillkor} it follows immediately that the group $G$ is invariant under the the Cartan involution
\begin{eqnarray}
\theta(g):=(g^*)^{-1} \label{involution}
\end{eqnarray}
and hence the space $G/K$ is equipped with the family of reflections $\{\sigma_{gK}\}_{gK \in G/K}$ given by
\begin{eqnarray}
\sigma_{gK}(xK):=g\theta(g^{-1}x)K
\end{eqnarray}
which furnish $G/K$ with the structure of a Riemannian symmetric space of the noncompact type.
In particular, for any $h \in B_1 (\mathbb{H})$, there is a unique geodesic joining 0 and $h$.
We let $\varphi_h$ denote the reflection in the midpoint, $m_h$, of this geodesic. The isometry $\varphi_h \in G$ is uniquely
characterised by the properties
\begin{eqnarray}
\varphi_h(m_h)&=&h,\\
d\varphi_h(m_h)&=&-Id_{T_{m_h}(B_1 (\mathbb{H}))}.
\end{eqnarray}
We let $Sp(1)_1$ and $Sp(1)_2$ denote the ``upper'' and ``lower'' subgroups of $K$ given by
\begin{eqnarray*}
Sp(1)_1&=&\left\{ \left( \begin{array}{cc} 
a & 0\\
0 & 1\\
\end{array}\right) \in K \right\},\\
Sp(1)_2&=&\left\{ \left( \begin{array}{cc} 
1 & 0\\
0 & d\\
\end{array}\right) \in K \right\}.
\end{eqnarray*}
For $k=\left( \begin{array}{cc} 
a & 0\\
0 & d\\
\end{array}\right) \in K$,
we will write $k=(a,d):=(k_1,k_2)$.
The group $Sp(1)_1 \cong SU(2)$ is then a normal subgroup of $K$. We will write
\begin{eqnarray}
\pi: K \rightarrow K/Sp(1)_1\cong SU(2)
\end{eqnarray}
for the natural projection onto the quotient group.

The group $K$ acts on the tangent space $T_0(B_1(\mathbb{H}))$ by the differentials at $0$ of the actions on $B_1(\mathbb{H})$.
By the restriction to the subgroup $Sp(1)_1$ we have a representation of $SU(2)$ on $T_0(B_1(\mathbb{H}))$.
We can define an $SU(2)$-representation, $\mu_h$, on the tangent space $T_h(B_1 (\mathbb{H}))$ for any $h$ by the formula
\begin{eqnarray}
\mu_h(l)v :=d\varphi_h(0) \circ dl(0) \circ d\varphi_h^{-1}(h)v, \label{SU(2)-str}\\
v \in T_h(B_1 (\mathbb{H})), l \in Sp(1)_1. \nonumber
\end{eqnarray}
The family $\{\mu_h\}_{h \in B_1(\mathbb{H})}$ of $SU(2)$-representations amounts to an action of $SU(2)$ as gauge transformations of
the tangent bundle $T(B_1 (\mathbb{H}))$. It is, however, not invariant under the action of $G$ as automorphisms of the bundle.
Indeed, if we define,for $h \in B_1(\mathbb{H}),  g \in G$, 
\begin{eqnarray}
\kappa_{g,h}:=\varphi_{g(h)}^{-1}g\varphi_h \in K,
\end{eqnarray}
then
\begin{eqnarray}
\mu_{g(h)}(l)dg(h)v=dg(h)\mu_h(\kappa_{g,h}^{-1}l\kappa_{g,h})v,
\end{eqnarray}
where the element $\kappa_{g,h}^{-1}l\kappa_{g,h}$ belongs to the subgroup $Sp(1)_1$ since it is normalised by $K$.
Hence the principal fibre bundle over $B_1(\mathbb{H})$  defined by
the family $\{\mu_h\}_{h \in B_1(\mathbb{H})}$ is $G$-equivariant, though not elementwise.
This shows that the symmetric space has a quaternionic structure and is a quaternionic symmetric space in the sense defined by 
Wolf (cf. \cite{wolf}). 

\subsection{Harish-Chandra realisation}
We consider again the embedding $\iota$ defined in eq. \eqref{matrisjota}. If we set 
\begin{eqnarray*}
G'&=&SU(2,2),\\
K'&=&S(U(2) \times U(2))\\
&:=&\left\{ \left( \begin{array}{cc}
A & 0\\
0 & D \\
\end{array}\right) \in SU(2,2)| A \in U(2), D \in U(2), \det(A)\det(D)=1\right\},
\end{eqnarray*}
$\iota$ induces an embedding of pairs $(G,K) \hookrightarrow (G',K')$ and hence descends to an embedding
\begin{eqnarray}
G/K \hookrightarrow G'/K'
\end{eqnarray}
of the corresponding symmetric spaces.
We will write $SU(2)_1$ and $SU(2)_2$ for the images $\iota(Sp(1)_1)$ and $\iota(Sp(1)_2)$ respectively.

The Hermitian symmetric 
space $G'/K'$ is by the Harish-Chandra realisation holomorphically, and $G$-equivariantly, equivalent to the bounded symmetric 
domain of type $I$
\begin{eqnarray}
G'/K' \cong \mathscr{D}:=\{Z \in M_2(\mathbb{C})|I_2-Z^*Z>0\}.
\end{eqnarray}
The action of $G'$ on $\mathscr{D}$ is given by
\begin{eqnarray}
g(Z)=(AZ+B)(CZ+D)^{-1}, \label{fraktionell G}
\end{eqnarray}
if $g=\left( \begin{array}{cc}
A & B \\
C & D \\
\end{array}\right)$ is a block matrix with blocks of size $2\times2$.
The symmetric space $G/K$ is thus embedded into $\mathscr{D}$ as the subset
\begin{eqnarray*}
\mathcal{D} :=\left\{Z \in M_2(\mathbb{C})|I_2-Z^*Z>0, Z=\left(\begin{array}{cc}
z & w \\
-\overline{w} & \overline{z} \\
\end{array}\right)\, z,w \in \mathbb{C}\right\},
\end{eqnarray*}
and the action is given by
\begin{eqnarray}
\iota(g)(\iota(h))=(\iota(a)\iota(h)+\iota(b))(\iota(c)\iota(h)+\iota(d))^{-1}=\iota(g(h)),
\end{eqnarray}
where $g(h)$ is the action defined in \eqref{fraktionell}.

For any $Z \in \mathscr{D}$, the tangent space $T_{Z}(\mathscr{D})$ is identified with the complex vector space $M_2(\mathbb{C})$ and
the differentials at $0$ of the $K'$ actions are given by
\begin{eqnarray}
dk(0)Z=AZD^{-1}, \qquad Z \in M_2(\mathbb{C}), k=\left( \begin{array}{cc}
A & 0 \\
0 & D \\
\end{array}\right) \in K'.
\end{eqnarray}

\subsection{Cartan subalgebra and root system}
Recall the Cartan involution $\theta$ on $G$ \eqref{involution}.
Its differential at the identity determines a decomposition of 
$\mathfrak{g}$ into the $\pm 1$-eigenspaces $\mathfrak{k}$ and $\mathfrak{p}$ respectively,
\begin{eqnarray}
\mathfrak{g}=\mathfrak{k} \oplus \mathfrak{p}, \label{cartanuppdelning}
\end{eqnarray}
where
\begin{eqnarray*}
\mathfrak{k}&=&\left\{\left(
\begin{array}{cc}
X & 0\\
0 & Y\\
\end{array}\right)| X, Y\in \mathbb{H}, X^*=-X, Y^*=-Y, \mbox{tr}X+\mbox{tr}Y=0\right\},\\
\mathfrak{p}&=&\left\{\left(
\begin{array}{cc}
0 & X\\
X^* & 0\\
\end{array}\right)| X \in \mathbb{H} \right\}.
\end{eqnarray*}
Let $\mathfrak{t} \subset \mathfrak{k}$ denote the subalgebra (realised as complex matrices)
\begin{eqnarray}
\mathfrak{t}=\left\{\left( \begin{array}{cccc}
si & 0 & 0 & 0\\
0 & -si & 0 & 0\\
0 & 0 & ti & 0\\
0 & 0 & 0 & -ti \\
\end{array} \right) | s, t \in \mathbb{R} \right\}.
\end{eqnarray}
It has a basis $\{H_1, H_2\}$, where
\begin{eqnarray}
H_1&=&\left( \begin{array}{cccc}
i & 0 & 0 & 0\\
0 & -i & 0 & 0\\
0 & 0 & 0 & 0\\
0 & 0 & 0 & 0 \\
\end{array} \right),\\
H_2&=&\left( \begin{array}{cccc}
0 & 0 & 0 & 0\\
0 & 0 & 0 & 0\\
0 & 0 & i & 0\\
0 & 0 & 0 & -i \\
\end{array} \right).
\end{eqnarray}
Let $\mathfrak{g}^{\mathbb{C}}$ be the complexification of $\mathfrak{g}$, and $\mathfrak{k}^{\mathbb{C}}$ and 
$\mathfrak{p}^{\mathbb{C}}$ denote the complexifications of $\mathfrak{k}$ and $\mathfrak{p}$ respectively.
The Cartan decomposition induces the decomposition
\begin{eqnarray}
\mathfrak{g}^{\mathbb{C}}=\mathfrak{k}^{\mathbb{C}} \oplus \mathfrak{p}^{\mathbb{C}}.
\end{eqnarray}
The complexification $\mathfrak{t}^{\mathbb{C}}\subset \mathfrak{k}^{\mathbb{C}}$ is a compact Cartan subalgebra of 
$\mathfrak{g}^{\mathbb{C}}$. 
Let $\Delta$ denote the set of roots, and for $\alpha \in \Delta$, we let $\mathfrak{g}^{\alpha}$ denote the corresponding root space.
Then, for each $\alpha \in \Delta$ either the inclusion
$\mathfrak{g}^{\alpha} \subseteq \mathfrak{k}^{\mathbb{C}}$ or the inclusion 
$\mathfrak{g}^{\alpha} \subseteq \mathfrak{p}^{\mathbb{C}}$ holds. In the first case, we call the root compact, and in the
second case we call it non-compact. Let $\Delta_{\mathfrak{k}}$, and $\Delta_{\mathfrak{p}}$ denote the set of compact roots and the
set of non-compact roots respectively.
We order the roots by letting the ordered basis $\{-\sqrt{-1}H_1^*, -\sqrt{-1}H_2^*\}$ for the real vector space 
$\sqrt{-1}\mathfrak{t}^*$ define a lexicographic ordering.
We let $\Delta_{\mathfrak{k}}^+$ denote the set of positive compact roots, and we let $\Delta_{\mathfrak{p}}^+$ denote the
set of positive non-compact roots.

The roots are given by
\begin{eqnarray}
\Delta_{\mathfrak{k}}&=&\{\pm 2\sqrt{-1}H_1^*, \pm2\sqrt{-1}H_2^*\},\\
\Delta_{\mathfrak{p}}&=&\{\pm \sqrt{-1}(H_1^*+H_2^*), \pm \sqrt{-1}(H_1^*-H_2^*)\}.
\end{eqnarray}
In terms of quaternionic matrices, the corresponding root spaces are
\begin{eqnarray}
\mathfrak{g}_{\pm 2\sqrt{-1}H_1^*}=\mathbb{C}  \left( \begin{array}{cc}
j \mp \sqrt{-1}k & 0\\
0 & 0\\ 
\end{array}\right),\\
\mathfrak{g}_{\pm 2\sqrt{-1}H_2^*}=\mathbb{C}\left( \begin{array}{cc}
0 & 0\\
0 & j \mp \sqrt{-1}k\\ 
\end{array}\right),
\end{eqnarray}
and 
\begin{eqnarray}
\mathfrak{g}_{\pm \sqrt{-1}(H_1^*+H_2^*)}&=&\mathbb{C}  \left( \left( \begin{array}{cc}
0 & j\\
-j & 0\\ 
\end{array}\right) \mp
\left( \begin{array}{cc}
0 & k\\
-k & 0\\ 
\end{array}\right) \right),\\
\mathfrak{g}_{\pm \sqrt{-1}(H_1^*-H_2^*)}&=&\mathbb{C}  \left( \left( \begin{array}{cc}
0 & 1\\
1 & 0\\ 
\end{array}\right) \mp
\left( \begin{array}{cc}
0 & i\\
-i & 0\\ 
\end{array}\right) \right)
\end{eqnarray}
respectively.

According to the lexicographic ordering on $\sqrt{-1}\mathfrak{t}^*$ determined by the ordered basis $\{-\sqrt{-1}H_1^*, -\sqrt{-1}H_2^*\}$,
the positive noncompact roots are 
\begin{eqnarray}
\alpha_1=-\sqrt{-1}H_1^*+\sqrt{-1}H_2^*,\\ 
\alpha_2=-\sqrt{-1}H_1^*-\sqrt{-1}H_2^*,
\end{eqnarray}
and $\alpha_1 < \alpha_2$.
Moreover, $\alpha_1+\alpha_2=-2\sqrt{-1}H_1^*$, i.e., the sum is a root. Hence $\{\alpha_1\}$ is a maximal sequence of strongly orthogonal 
positive noncompact roots.
We let $B(\cdot, \cdot)$ denote the Killing form on $\mathfrak{g}$. We use it to identify $\mathfrak{g}^*$ with $\mathfrak{g}^*$
according to
\begin{eqnarray}
\alpha(X):=B(X, H_{\alpha}), \alpha \in \mathfrak{g}^*, X \in \mathfrak{g}^*.
\end{eqnarray}
Via this identification, the Killing form induces a bilinear form on $\mathfrak{g}^*$ by
\begin{eqnarray}
\langle \alpha, \beta \rangle:=B(H_{\alpha}, H_{\beta}).
\end{eqnarray}
For $\alpha \in \Delta$, we select a root vector $E_{\alpha} \in \mathfrak{g}_{\alpha}$ in such a way that 
\begin{eqnarray}
B(E_{\alpha}, E_{-\alpha})=\frac{2}{\langle \alpha, \alpha \rangle}.
\end{eqnarray}

\subsection{Iwasawa decomposition}
Consider the maximal abelian subspace 
\begin{eqnarray}
\mathfrak{a}:=\mathbb{R}(E_{\alpha_1}+E_{-\alpha_1})=\left\{\left(
\begin{array}{cc}
0 & tI_2\\
tI_2 & 0\\ \end{array}\right)|\, t \in \mathbb{R}\right\}
\end{eqnarray}
of $\mathfrak{p}$.
The Iwasawa decomposition of $\mathfrak{g}$ with respect to $\mathfrak{a}$ is given by
$$\mathfrak{g}=\mathfrak{n} \oplus \mathfrak{a} \oplus \mathfrak{k}.$$
The corresponding global decomposition is $$G=NAK,$$
where, written as quaternionic matrices
\begin{eqnarray*}
N&=&\left\{\left(
\begin{array}{cc}
1+q & -q\\
q & 1-q\\
\end{array}\right)|\, q \in \mathbb{H}, q^*=-q \right\},\\
A&=&\left\{\left(
\begin{array}{cc}
\cosh t & \sinh t\\
\sinh t & \cosh t\\
\end{array}\right)| \,t \in \mathbb{R} \right\}.
\end{eqnarray*}

\begin{remark}
One can just as well use an Iwasawa decomposition $G=KAN$, the correpondence between these two decompositions being
$(nak)^{-1}=k^{-1}a^{-1}n^{-1}$. In the sequel will see that it is sometimes convenient use this other decomposition 
as a means for finding the components in our decomposition.
\end{remark}
In the sequel we will need the explicit formulas for the $NAK$-factorisation
\begin{eqnarray}
g=n(g)a(g)\kappa(g)
\end{eqnarray}
of an element $g \in Sp(1,1)$.
\begin{lemma} \label{faktorer}
For $g=\left(
\begin{array}{cc}
a & b\\
c & d\\
\end{array}\right)$, and $\log a(g)=t(E_{\alpha_1}+E_{-\alpha_1})$,  $e^t$ and $\kappa(g)$ are given by
\begin{eqnarray*}
e^{t}&=&
\frac{(1-|bd^{-1}|^2)^{1/2}}{|1-bd^{-1}|},\\
\kappa(g)&=&
e^{t}\left(
\begin{array}{cc}
a-c & 0\\
0 & d-b\\
\end{array}\right).
\end{eqnarray*}
\end{lemma}

\begin{proof}
The proof is by straightforward computation. We prove only the second statement.

The identity 
$$\left(
\begin{array}{cc}
a & b\\
c & d\\
\end{array}\right)=
\left(
\begin{array}{cc}
1+q & -q\\
q & 1-q\\
\end{array}\right)
\left(
\begin{array}{cc}
\cosh t & \sinh t\\
\sinh t & \cosh t\\
\end{array}\right)
\left(
\begin{array}{cc}
u_1 & 0\\
0 & u_2\\
\end{array}\right)$$
is equivalent to
\begin{eqnarray*}
&&\left(\begin{array}{cc}
a & b\\
c & d\\
\end{array}\right)\\
&&=\left(\begin{array}{cc}
(\cosh t+q(\cosh t-\sinh t))u_1 &(\sinh t+q(\sinh t-\cosh t))u_2 \\
(\sinh t+q(\cosh t-\sinh t))u_1&(\cosh t+q(\sinh t-\cosh t))u_2 \\
\end{array}\right).
\end{eqnarray*}
Hence 
\begin{eqnarray*}
a-c&=&e^{-t}u_1\\
d-b&=&e^{-t}u_2.
\end{eqnarray*}

\end{proof}

\section{Quaternionic discrete series representations}

\subsection{Generalities}
The Cartan decomposition \eqref{cartanuppdelning}
decomposes $\mathfrak{g}$ into two invariant subspaces for the adjoint action of $K$.
Moreover, we have the isomorphism of $K$-representations
\begin{eqnarray}
\mbox{Ad}_{\mathfrak{g}/\mathfrak{k}} \cong \mbox{Ad}_{\mathfrak{p}}.
\end{eqnarray}
We extend $\mbox{Ad}_{\mathfrak{p}}$ to a complex linear representation of $K$ on 
the space $(\mathfrak{g}/\mathfrak{k})^{\mathbb{C}} \cong \mathfrak{p}^{\mathbb{C}}$.

Consider now the surjective mapping
\begin{eqnarray}
p: G \rightarrow G/K, p(g)=gK.
\end{eqnarray}
The differential at the origin 
\begin{eqnarray}
dp(e): \mathfrak{g} \rightarrow T_{eK}(G/K) \label{diffp}
\end{eqnarray}
intertwines the adjoint action of $K$ on $\mathfrak{g}$ with the differential action on the tangent space $T_{eK}(G/K)$.
The kernel of $dp(e)$ is $\mathfrak{k}$, and as $K$-representations we thus have the isomorphism
\begin{eqnarray}
\mbox{Ad}_{\mathfrak{p}}^* \cong (dK(o)^{\mathbb{C}})^*,
\end{eqnarray}
where the right hand side denotes the complex linear dual to the representation given by the complexified actions of the 
tangent maps at the origin.
Using the quotient mapping induced by \eqref{diffp} and the realisation of the differential action of $K$ at the tangent space 
$T_0(\mathscr{D})$, we obtain the formula
\begin{eqnarray}
\mbox{Ad}_{\mathfrak{p}}^*(k)Z=(A^{-1})^tZD^t.
\end{eqnarray}
Here $Z \in M_2(\mathbb{C}) \cong T_0^*(\mathscr{D}) \cong (T_0(\mathcal{D})^{\mathbb{C}})^*$, 
and $k=\left(\begin{array}{cc}
A & 0\\
0 & D\\
\end{array}\right) \in K \cong SU(2) \times SU(2)$.
The restriction, $\mbox{Ad}_{\mathfrak{p}}^*|_{SU(2)_1}$, to the subgroup $SU(2)_1$ is then given by
\begin{eqnarray}
\mbox{Ad}_{\mathfrak{p}}^{*}|_{SU(2)_1}(k)Z=(A^{-1})^tZ.
\end{eqnarray}
If we let $\{E_{ij}\}$, $i,j=1,2$ denote the standard basis for the complex vector space (i.e., $E_{ij}$ has $1$ at the position
on the $i$th row and $j$th column and zeros elsewhere), then clearly the subspace
\begin{eqnarray}
V:=\mathbb{C}E_{11} \oplus \mathbb{C}E_{21} \cong \mathbb{C}^2
\end{eqnarray}
spanned by the basis elements in the first column is $SU(2)_1$-invariant. Likewise, the subspace spanned by the basis elements
of the second column is invariant.
We now let $\tau$ denote the representation given by restricting the $K$-representation $\mbox{Ad}_{\mathfrak{p}}^*|_{SU(2)_1}$ to the subspace $V$, 
and let $\tau_k$ denote the $k$th symmetric tensor power of the representation $\tau$.
Then clearly, the natural identification of $\tau_k$ with a representation of $SU(2)$ is equivalent to the standard representation of
$SU(2)$ on the space of polynomial functions $p(z,w)$ on $\mathbb{C}^2$ of homogeneous degree $k$, i.e., we have 
\begin{eqnarray}
\tau_k(l_1,l_2)p(z,w):=p(l_1^{-1}(z,w)):=p(az+\overline{b}w,-bz+\overline{a}w), \label{polynomrum}
\end{eqnarray}
where $l_1^{-1}=\left( \begin{array}{cc}
a & \overline{b}\\
-b & \overline{a}
\end{array}\right) \in SU(2)$.
We let $V_{\tau_k}$ denote the representation space for $\tau_k$.
The (smoothly) induced representation $\mbox{Ind}_K^G(\tau_k)$ is then defined on the space 
\begin{eqnarray*}
C^{\infty}(G, \tau_k):=\{f \in C^{\infty}(G, V_{\tau_k})|f(gl^{-1})=\tau_k(l)f(g) \,\,\forall g \in G \, \forall l \in K\}, \label{indfunkt}
\end{eqnarray*}
i.e., on the space of smooth sections on the $G$-homogeneous vector bundle
\begin{eqnarray}
\mathcal{V}^k \rightarrow G/K:=G \times_K V_{\tau_k} \rightarrow G/K.
\end{eqnarray}
We fix the $K$-invariant inner product on $ \langle , \rangle_k$ on $V_{\tau_k}$ given by
\begin{eqnarray}
\langle p, q \rangle_k:=p(\partial)(q^*)(0),
\end{eqnarray}
where $p(\partial)$ is the differential operator defined by substituting $\frac{\partial}{\partial z}$ for $z$, and 
$\frac{\partial}{\partial w}$ for $w$ in the polynomial function $p(z,w)$, and 
\begin{eqnarray*}
(\sum_{j=1}^ka_jz^jw^{k-j})^*:=\sum_{j=1}^k \overline{a_j}z^jw^{k-j}).
\end{eqnarray*}
We use this inner product to define an Hermitian metric on $\mathcal{V}^k$ by 
\begin{eqnarray} 
h_Z(u,v):=\langle (g^{-1})_Zu,(g^{-1})_Zv \rangle_k, \qquad u,v \in \mathcal{V}_Z^k,
\end{eqnarray}
where $Z=gK$ and $(g^{-1})_Z$ denotes the fibre map $\mathcal{V}_Z^k \rightarrow \mathcal{V}_0^k \cong V_k$ associated with $g^{-1}$.
For a fixed choice, $\iota$, of $G$-invariant measure on $G/K$ we define $L^2(\mbox{Ind}_K^G(\tau_k))$ as the Hilbert space
completion of the space 
\begin{eqnarray}
\left \{s \in \Gamma(G/K, \mathcal{V}^k)| \int_{G/K}h_Z(s,s)d\iota(Z) < \infty \right\}. \label{sektionsnorm}
\end{eqnarray}
The tensor product representation $\tau_k \otimes \mbox{Ad}(K)|_{\mathfrak{p}^{\mathbb{C}}}$ decomposes into $K$-types according to
\begin{eqnarray}
\tau_k \otimes \mbox{Ad}(K)|_{\mathfrak{p}^{\mathbb{C}}}=\sum_{\beta \in \Delta_{\mathfrak{p}}}m_{\beta}\pi_{\beta-k\sqrt{-1}H_1^*},
\end{eqnarray}
where $m_{\beta} \in \{0,1\}$, and $\pi_{\beta-k\sqrt{-1}H_1^*}$ is the irreducible representation of $K$ with highest weight 
$\beta-k\sqrt{-1}H_1^*$.
Let $\tau_k^-$ be the subrepresentation of the tensor product given by 
\begin{eqnarray}
\tau_k^-=\sum_{\beta \in \Delta_{\mathfrak{p}}^-}m_{\beta}\pi_{\beta-k\sqrt{-1}H_1^*},
\end{eqnarray}
and let $V_k^-$ be the subspace of $V_{\tau_k} \otimes \mathfrak{p}^{\mathbb{C}}$ on which $\tau_k^-$ operates.
Let $P: V_{\tau_k} \otimes \mathfrak{p}^{\mathbb{C}} \rightarrow V_k^-$ be the orthogonal projection.
Define the space $C^{\infty}(G, \tau_k^-)$ in analogy with \eqref{indfunkt}. 
We recall that the \emph{Schmid  $D$ operator} is a differential operator mapping the space $C^{\infty}(G, \tau_k)$ into 
$C^{\infty}(G, \tau_k^-$) and is defined as
\begin{eqnarray}
Df(g)=\sum_i P(X_if(g) \otimes X_i),
\end{eqnarray}
where $\{X_i\}$ is any orthonormal basis for $\mathfrak{p}^{\mathbb{C}}$, and $X_if$ denotes left invariant differentiation, i.e., 
\begin{eqnarray*}
Xf(g)&:=&\frac{d}{dt}f(g \exp(tX))|_{t=0}, X \in \mathfrak{p},\\
Zf(g)&:=&Xf(g)+iY(g), Z=X+iY \in \mathfrak{p}^{\mathbb{C}}.
\end{eqnarray*}
The subspace $\ker D \cap L^2(\mbox{Ind}_K^G(\tau_k))$ is then invariant under the left action of $G$ and defines
an irreducible representation of $G$ belonging to the quaternionic discrete series. We let $\mathcal{H}_k$ denote this
representation space. By \cite{GW}, it belongs to the discrete series for $k \geq 1$.

\begin{remark}
The model we use to describe the Hilbert space $\mathcal{H}_k$ can be used to realise any discrete series representation by
induction from $K$ to $G$ of the minimal $K$-type for any pair $(G,K)$ where $G$ is semisimple and $K$ is maximal compact (cf. \cite{knapp1}). 
By \cite{GW}, for $k \geq 1$, $\tau_k$ occurs as a minimal $K$-type for some
discrete series representation of $G=Sp(1,1)$.
\end{remark}

\subsection{Global trivialisation}
Let us for a while view the representation space $\mathcal{H}_k$ as a space of sections of the vector bundle $\mathcal{V}^k \rightarrow G/K$.
We recall the diffeomorphism $G/K \cong \mathcal{D}$ given by $gK \mapsto g \cdot 0$. This lifts to a global trivialisation, $\Phi$, of
the bundle $\mathcal{V}^k \rightarrow G/K$ given by
\begin{eqnarray}
\qquad \Phi: G \times_K V^{\tau_k} \rightarrow \mathcal{D} \times V_{\tau_k},\,\,
\Phi([(g,v)]):= (g \cdot 0, \tau_k(J(g,0))v), \label{globaltrivial}
\end{eqnarray}
where $J(g,Z)$ denotes the $K^{\mathbb{C}}$-component of $g\exp Z$ - the \emph{automorphic factor of g at Z} (cf. \cite{satake}).

If $F : G \rightarrow V^{\tau_k}$ is a function in $C^{\infty}(G, \tau_k)$, its trivialised counterpart is the function 
$f: \mathcal{D} \rightarrow V_{\tau_k}$ given by
\begin{eqnarray}
f(g\cdot 0):=\tau_k(J(g,0))F(g).
\end{eqnarray}
In the trivialised picture, the group $G$ acts on functions on $\mathcal{D}$ by
\begin{eqnarray}
gf(Z):=\tau_k(J(g^{-1},Z))^{-1}f(g^{-1}Z).
\end{eqnarray}
More explicitly, if $g^{-1}=\left(\begin{array}{cc}
A & B\\
C & D\\
\end{array}\right)$ (considered as a matrix in $SU(2,2)$), then
\begin{eqnarray*}
J(g^{-1}Z)=\left( \begin{array}{cc}
A-(AZ+B)(CZ+D)^{-1}C & 0\\
0 & D\\
\end{array}\right) \in SL(4,\mathbb{C}),
\end{eqnarray*}
and 
\begin{eqnarray*}
gf(Z)=\stackrel{k}{\odot} 
(A-(AZ+B)(CZ+D)^{-1}C)f((AZ+B)(CZ+D)^{-1}).
\end{eqnarray*}
The action of $SU(2)$ on the vector space $V_{\tau_k}$ is here naturally extended to an action of $SL(2,\mathbb{C})$ by the formula
\eqref{polynomrum}.

In the trivialised picture, the norm \eqref{sektionsnorm} can also be described explicitly.
\begin{prop} Let $k \geq 1$.
In the realisation of the Hilbert space $\mathcal{H}_k$ as a space of $V_{\tau_k}$-valued functions on $\mathcal{D}$, the norm \eqref{sektionsnorm}
is given by
\begin{eqnarray}
\|f\|_k:=\int_{B_1(\mathbb{H})}(1-|q|^2)^k\langle f(q),f(q)\rangle_k(1-|q|^2)^{-4}dm(q).
\end{eqnarray}
\end{prop}

\begin{proof}
For $Z=gK$, a fibre map $(g^{-1})_Z: V_{\tau_k} \rightarrow V_{\tau_k}$ is given by
\begin{eqnarray}
(g^{-1})_Zv=\tau_k(J(g,0))^{-1}v.
\end{eqnarray}
If $g=\left(\begin{array}{cc}
\cosh t I_2 & \sinh t I_2\\
\sinh t I_2 & \cosh t I_2 \\
\end{array}\right)$, the automorphic factor $J(g,0)$ is given by (cf.\cite{knapp1} )
\begin{eqnarray}
&&J(g,0) \label{realautomorf} \\
&&=\left(\begin{array}{cc}
\cosh t ^{-1}I_2 & 0\\
0 & \cosh t I_2\\
\end{array}\right)
=\left(\begin{array}{cc}
(1-\tanh^2 t)^{1/2} I_2 & 0\\
0 & \cosh t I_2 \\
\end{array}\right). \nonumber
\end{eqnarray}
A general point $Z \in \mathcal{D}$ can be described as $Z=kgK$ for $g$ as above. The cocycle condition
\begin{eqnarray}
J(kg,0)=J(k,g0)J(g,0)
\end{eqnarray}
then implies that
\begin{eqnarray}
J(kg,0)=\left(\begin{array}{cc}
k_1(1-\tanh ^2 t)^{1/2}I_2 & 0\\
0 & k_2\cosh t I_2 \\
\end{array}\right),
\end{eqnarray}
if $k=(k_1,k_2) \in SU(2) \times SU(2)$.
Hence, for $k=1$
\begin{eqnarray*}
h_Z(u,v)&=&\mbox{tr} \left((k_1(1-\tanh^2 t)^{1/2}I_2)^t u((k_1(1-\tanh^2 t)^{1/2}I_2)^tv)^*\right)\\
&=& \mbox{tr}\left((I_2-ZZ^*)^tuv^*\right).
\end{eqnarray*}
For arbitrary $k$, we have
\begin{eqnarray}
h_Z(u,v)&=&\mbox{tr}\left( \stackrel{k}{\odot}(I_2-ZZ^*)^tuv^*\right).
\end{eqnarray}
By analogous considerations, it follows that the invariant measure is given by
\begin{eqnarray*}
d\iota(Z)=\det(I_2-Z^*Z)^{-2}dm(Z),
\end{eqnarray*}
where $dm(Z)$ denotes the Lebesgue measure.
Hence, we obtain the formula
\begin{eqnarray}
\int_{\mathcal{D}}\langle \stackrel{k}{\odot}(I_2-ZZ^*)^tf(Z),f(Z) \rangle_k \det(I_2-Z^*Z)^{-2}dm(Z)
\end{eqnarray}
for the norm \eqref{sektionsnorm}. In quaternionic notation, this translates into the statement of the proposition.
\end{proof}


\section{Principal series representations and the Szegö map}
In this section we will consider a realisation of the discrete series representation $L^2(\mbox{Ind}K^G(\tau_k))$
as a quotient of a certain nonunitary principal series representation. We first state the theorem, and then
we investigate how the given principal series representation decomposes into $K$-types.
From now on we fix the number $k$ and simply write $\tau$ for $\tau_k$.

Recall the maximal abelian subspace $\mathfrak{a}$ of $\mathfrak{p}$ and consider the parabolic subgroup 
$$P=MAN$$ of $G$, where
\begin{eqnarray*}
M=Z_{K}(A)=\left\{\left(
\begin{array}{cc}
u & 0\\
0 & u\\
\end{array}\right)|\, u \in SU(2)\right\},
\end{eqnarray*}
and $A$ and $N$ are the ones that occur in the Iwasawa decomposition.
Let $\sigma$ be the restriction of the representation $\tau$ to the subgroup $M$.
Then, clearly, the subspace defined by the $M$-span of the $\tau$-highest weight-vector equals $V_{\tau}$ and the
representation $\sigma$ is also irreducible. We will hereafter denote this representation space by $V^{\sigma}$.
Recall the identification of $V^{\sigma}$ with a space of homogeneous polynomials. We thus adopt a somewhat 
abusive notation and write $z^{\sigma}$ for the
highest weight-vector.
Let $\nu \in \mathfrak{a}^*$ be a real-valued linear functional and consider the representation 
\begin{eqnarray}
\sigma \otimes \exp (\nu) \otimes 1
\end{eqnarray}
of $P$.
The induced representation $\mbox{Ind}_P^G(\sigma \otimes \exp (\nu) \otimes 1)$ is defined on the set
of continuous functions $f: G \rightarrow V_{\sigma}$ having the $P$-equivariant property
\begin{eqnarray}
f(gman)=e^{-\nu(\log a)}\sigma(m)^{-1}f(g).\label{eqindP}
\end{eqnarray}
The action of $G$ on this space is given by left translation,
$$\mbox{Ind}_P^G(\sigma \otimes \exp (\nu) \otimes 1)(f)(x):=L_{g^{-1}}f(x)
=f(g^{-1}x).$$
Consider now the smoothly induced representation $\mbox{Ind}_{M}^{K}(\sigma)$ which operates on the space, $C^{\infty}(K,\sigma)$, 
of all smooth functions
$f: K \rightarrow V_{\sigma}$ having the $M$-equivariance property 
\begin{eqnarray}
f(km)=\sigma(m)^{-1}f(k)\label{eqind} 
\end{eqnarray}
with $K$-action given by left translation.
The Iwasawa decomposition $G=KAN$ shows that, a fortiori, $G=KMAN$ (although this factorisation is not unique).
Given a linear functional $\nu \in \mathfrak{a}^*$, we can therefore extend any such function on $K$ to a 
function on $G$ by setting 
$$f(kman)=e^{-\nu(\log a)}\sigma(m)^{-1}f(k), \mbox{for}\,\, g=kman.$$
The equivariance property \eqref{eqind} of $f$ guarantees that this is indeed well-defined even though the factorisation
of $g$ is not. The extended function $f$ has the $P$-equivariance property \eqref{eqindP}. In fact, this extension
procedure defines a bijection between the representation spaces of the representations
$\mbox{Ind}_{M}^{K}(\sigma)$ and $\mbox{Ind}_P^G(\sigma \otimes \exp (\nu) \otimes 1)$.
There is a natural pre-Hilbert space structure on this representation space given by
$$\|f\|^2=\int_K\|f(k)\|_{\sigma}^2dk,$$ where 
$\|\cdot\|_{\sigma}$ denotes the inner product on $V_{\sigma}$ and $dk$ is the Haar measure on $K$.
The completion of the space of $M$-equivariant smooth functions $K \rightarrow V_{\sigma}$ with respect to
this sesquilinear form can be identified with the space of all square-integrable $V_{\sigma}$- valued functions having the
property \eqref{eqind}. We will denote the $K$-representation on this space by $L^2(\mbox{Ind}_{M}^{K}(\sigma))$.
By the extension procedure using $\nu$ described above, this completion can be extended to the space of all
$P$-equivariant $V_{\sigma}$-valued functions on $G$ such that the restriction to $K$ is square-integrable.

We now state the theorem by Knapp and Wallach.
\begin{thm}[\cite{K-W}, Thm. 6.1]
The \emph{Szegö mapping} with parameters $\tau$ and $\nu$ given by
\begin{eqnarray}
S(f)(x):=\int_K e^{\nu \log a(lx)}\tau(\kappa(lx)^{-1})f(l^{-1})dl
\end{eqnarray}
carries the the space $C^{\infty}(K, \sigma)$ into $C^{\infty}(G, \tau) \cap \ker D$, provided that
$\nu$ and $\tau$ are related by the formula
\begin{eqnarray}
\nu(E_{\alpha_1}+E_{-\alpha_1})=\frac{2 \langle -k\sqrt{-1}H_1^*+n_1\alpha_1, \alpha_1 \rangle}{\langle \alpha_1, \alpha_1 \rangle },
\label{parameterrelation}
\end{eqnarray}
where
\begin{eqnarray*}
n_1=|\{\gamma \in \Delta_{\mathfrak{p}}^+| \alpha(\gamma)=\alpha_1 \, \mbox{and}\,\, \alpha_1+\gamma \in \Delta\}|.
\end{eqnarray*}
\end{thm}
In this case $n_1=1$, since the root $\alpha_2$ is the only one satisfying the above condition.
Moreover, an easy calculation gives that
\begin{eqnarray}
E_{\alpha_1}=  \frac{1}{2}\left( \left( \begin{array}{cc}
0 & 1\\
1 & 0\\ 
\end{array}\right) +
\left( \begin{array}{cc}
0 & i\\
-i & 0\\ 
\end{array}\right) \right).
\end{eqnarray}
Hence, the condition \eqref{parameterrelation} takes the form
\begin{eqnarray}
\nu \left(\left( \begin{array}{cc}
0 & 1\\
1 & 0\\ 
\end{array}\right)\right)=k+2.
\end{eqnarray}
Hereafter, we will make the identification
\begin{eqnarray}
\nu=k+2 \label{funktionalidentifikation}
\end{eqnarray}
of the functional with a natural number.
We now proceed with a more detailed study of the representation $L^2(Ind_M^K(\sigma))$. 
\begin{lemma}\label{K-hom}
The representation $L^2(Ind_M^K(\sigma))$ is $K$-equivalent to \\
$L^2(K/M) \otimes V_{\sigma}$.
\end{lemma}

\begin{proof}
Let $f$ be a continuous function from $K$ to $V_{\sigma}$ having the property of $M$-equivariance
$$f(km)=\sigma(m)^{-1}f(k), k \in K, m \in M.$$
Then the function $$\tilde{f}(k):= \tau(k)f(k)$$ is clearly right $M$-invariant and hence we can define the function
$F: K/M \rightarrow V_{\sigma}$ by
$$F(kM)=\tilde{f}(k).$$
This is obviously well-defined. By choosing a basis $\{e_j\}$ for $V_\sigma$, we can write
$$F(kM)=\Sigma_j F_j(kM)e_j$$ for some complex-valued functions $F_j$. 
We now define a mapping $$T:L^2(Ind_M^K(\sigma)) \rightarrow L^2(K/M) \otimes V_{\sigma}$$ by
$$Tf:=\Sigma_j F_j \otimes e_j.$$
To see that this mapping is a bijection, note that any vector in the Hilbert space $L^2(K/M) \otimes V_{\sigma}$ can be 
uniquely expressed in the form $\Sigma_j G_j \otimes e_j.$ We can thus define a mapping
$$S: L^2(K/M) \otimes V_{\sigma} \rightarrow L^2(Ind_M^K(\sigma))$$ by
$$S(\Sigma_j G_j \otimes e_j)(k):=\tau(k)^{-1}\Sigma_j g_j(kM)e_j$$ and it is easy to see that $S$ is the inverse of $T$. 

It remains now only to prove the $K$-equivariance. Pick therefore any element\\
$\Sigma_j G_j \otimes e_j$ from the Hilbert space on the right hand side. We have
$$k(\Sigma_j G_j \otimes e_j)=\Sigma_j G_j \circ L_{k^{-1}} \otimes \sigma(k)e_j.$$
If we denote the matrix coefficients of $\sigma(k)$ with respect to the basis $\{e_j\}$ by $\sigma(k)_{ij}$, we
have $$\sigma(k)e_j=\Sigma_i\sigma(k)_{ij}e_i$$ and hence 
$$\Sigma_j G_j \circ L_{k^{-1}} \otimes \sigma(k)e_j=\Sigma_{i,j}G_j \circ L_{k^{-1}} \otimes \sigma(k)_{ij}e_i.$$
Applying $S$ to the above expression yields 
\begin{eqnarray*}
S(\Sigma_{i,j}G_j \circ L_{k^{-1}} \otimes \sigma(k)_{ij}e_i)(k')&=&
\sigma(k')^{-1}\Sigma_{i,j}G_j(k^{-1}k'M)\sigma(k)_{ij}e_i\\
&=& \sigma(k')^{-1}\sigma(k)\Sigma_{j}G_j(k^{-1}k'M)e_j\\
&=&S(\Sigma_j G_j \otimes e_j) \circ L_{k^{-1}} (k').
\end{eqnarray*}
\end{proof}
We shall now examine the left action of $K$ on the $L^2(K/M)$-factor in the tensor product more closely. In particular, we 
are interested in a certain $K$-invariant subspace defined by a subclass of the $K$-types occurring in $L^2(K/M)$.
We recall the identification of the $K$-representation $\tau_j$ with a standard representation of $SU(2)$. We therefore
let $\tau_j$ also denote the corresponding $SU(2)$-representation, and we let $V_j$ denote the associated vector space
of polynomials. Any irreducible representation of $K=SU(2) \times SU(2)$ is isomorphic to a tensor product of irreducible
$SU(2)$-representations, i.e., it is realised on a space
\begin{eqnarray}
V_j^* \otimes V_i,
\end{eqnarray}
for some $i, j \in \mathbb{N}$.
With the fixed ordering of the roots, the polynomial function $(z,w) \mapsto z^j$ is a highest weight vector in $V_j$, and
the polynomial function $(z,w) \mapsto w^j$ is a lowest weight vector.
We will use the abusive notation where they are denoted by $z^j$ and $w^j$ respectively.
\begin{prop}\label{algsum}
The algebraic sum 
\begin{eqnarray}
W:=\bigoplus_{j \in \mathbb{N}}V_j^* \otimes V_{\sigma+j}
\end{eqnarray}
of $K$-types is a subspace of $L^2(\mbox{Ind}_M^K(\sigma)$. The highest weight vector for the $K$-type $V_j^* \otimes V_{\sigma +j}$
is given by the function
\begin{eqnarray}
f_j(k):= \langle \tau_j \circ \pi(k)z^j, w^j \rangle_j \,\tau(k)^{-1}z^\sigma.
\end{eqnarray}
\end{prop}

\begin{proof}
For $$k=\left(
\begin{array}{cc}
u_1 & 0\\
0 & u_2 \\
\end{array}\right), k'=\left(
\begin{array}{cc}
u'_1 & 0\\
0 & u'_2 \\
\end{array}\right)
 \in K,$$
we have
\begin{eqnarray*}
k'kM= \left(
\begin{array}{cc}
u'_1u_1 & 0\\
0 & u'_2u_2 \\
\end{array}\right)
M&=&
\left(
\begin{array}{cc}
u'_1u_1u_2^{-1}(u'_2)^{-1} & 0\\
0 & I \\
\end{array}\right)
\left(
\begin{array}{cc}
u'_2u_2 & 0\\
0 & u'_2u_2 \\
\end{array}\right)M\\
&=&
\left(
\begin{array}{cc}
u'_1u_1u_2^{-1}(u'_2)^{-1} & 0\\
0 & I \\
\end{array}\right)
M
\end{eqnarray*}
and this shows that the left action of $K=SU(2) \times SU(2)$ on $L^2(SU(2))$-functions is equivalent to
the action $L_{g^{-1}} \otimes R_{h}$:
$$(L_{g^{-1}} \otimes R_h)(g,h)f(l):=f(g^{-1}lh)$$
Then, by the Peter-Weyl Theorem, $L^2(K/M)$ decomposes into $K$-types according to
\begin{eqnarray}
L^2(K/M)\simeq \bigoplus_{j \in \widehat{SU(2)}}(V_j \otimes V^*_j).
\end{eqnarray}
Tensoring with $V_{\sigma}$ gives the sequence of $K$-isomorphisms
\begin{eqnarray*}
L^2(K/M) \otimes V_{\sigma}&\simeq&\bigoplus_{j \in \widehat{SU(2)}}(V_j \otimes V^*_j)\otimes V_{\sigma}\\
&\simeq& \bigoplus_{j \in \widehat{SU(2)}}(V^*_j \otimes V_j)\otimes V_{\sigma}\\
&\simeq& \bigoplus_{j \in \widehat{SU(2)}}V^*_j \otimes (V_j\otimes V_{\sigma}).
\end{eqnarray*}
Moreover, each term $(V_j\otimes V_{\sigma})$ has a \emph{Clebsch-Gordan}-decomposition
$$(V_j\otimes V_{\sigma})\simeq (V_{\sigma + j} \oplus \cdots)$$ and therefore each term 
$V_j \otimes V_{\sigma+j}$ will constitute a $K$-type in $L^2(K/M) \otimes V_{\sigma}$.
Such a $K$-type has a highest weight-vector $(w^j)^* \otimes z^{\sigma+j}$. Using first the embedding into 
$V_j^* \otimes (V_j\otimes V_{\sigma})$ and then the $K$-isomorphism given by Lemma \ref{K-hom}, we see that
highest weight-vector corresponds to the $M$-equivariant function
\begin{eqnarray}
f_j(k):= \langle \tau_j \circ \pi(k)z^j, w^j \rangle_j \,\sigma(k)^{-1}z^\sigma. 
\end{eqnarray}
\end{proof}

\section{Realisation of K-types}
By \cite{GW}, the only $K$-types occurring in the quaternionic discrete series for $Sp(1,1)$ are the ones that 
form the subspace $W$ in Proposition \ref{algsum}. In this section we compute their realisations as $V_{\tau}$-valued functions 
on $B_1(\mathbb{H})$ 
when restricted
to the submanifold
\begin{eqnarray}
A \cdot 0=\{t \in \mathbb{H}| -1 < t <1\}
\end{eqnarray} \label{realdel}
of $B_1(\mathbb{H})$.
For $s \in \mathbb{R}$, we let
\begin{eqnarray}
a_s=\left( \begin{array}{cc}
\cosh s & \sinh s\\
\sinh s & \cosh s\\
\end{array}\right) \in Sp(1,1).
\end{eqnarray}
Then $a_s \cdot 0=\tanh s \in A \cdot 0$. We start by computing the Szegö images of the $f_j$ when restricted to points $a_s$. 

Each of the standard $SU(2)$-representations, $V_N$, can be naturally extended to a representation of $GL(2,\mathbb{C})$ by
\begin{eqnarray}
p \mapsto p \circ g^{-1}, \qquad p \in V_n, g \in GL(2,\mathbb{C}).
\end{eqnarray}
This action of $GL(2,\mathbb{C})$ will occur frequently in the sequel.

\begin{lemma} \label{Szegö}
The Szegö transform of the highest weight-vector $f_j$ is given by
\begin{eqnarray*}
&&Sf_j(a_s)=(\cosh s)^{-\nu}\\
&&\times \int_{SU(2)}(\det(1-l \tanh s ))^{-(\nu+\sigma)/2}
\langle \tau_j(l^{-1})z^j,w^j \rangle_j
\,\sigma(1-l\tanh s)z^{\sigma}dl
\end{eqnarray*}
when restricted to the $A$-component in the decomposition $G=NAK$.
\end{lemma}

\begin{proof}
Take
\begin{eqnarray*}
k&=&\left(\begin{array}{cc}
u_1 & 0\\
0 & u_2 \\
\end{array}\right) \,\,\mbox{and}\\
x&=& \left(
\begin{array}{cc}
\cosh s & \sinh s\\
\sinh s & \cosh s\\
\end{array}\right).
\end{eqnarray*}
Then $$kx=\left(
\begin{array}{cc}
u_1\cosh s & u_1\sinh s\\
u_2\sinh s & u_2\cosh s\\
\end{array}\right),$$
and Lemma \ref{faktorer} gives that
\begin{eqnarray*}
e^{\nu(\log H(kx))}&=&
\left(\frac{ 1-|u_1u_2^{-1}\tanh s|^2}{|1-u_1u_2^{-1}\tanh s|^2}\right)^{\nu/2},\\
\kappa(kx)&=&\left(\frac{ 1-|u_1u_2^{-1}\tanh s|^2}{|1-u_1u_2^{-1}\tanh s|^2}\right)^{1/2}\\
&&\times
\left(\begin{array}{cc}
u_1\cosh s-u_2\sinh s & 0\\
0 & u_2\cosh s-u_1\sinh s \\
\end{array}\right).
\end{eqnarray*}
Hence $$\tau(\kappa(lx))^{-1}=\left(\frac{ 1-|u_1u_2^{-1}\tanh s|^2}{|1-u_1u_2^{-1}\tanh s|^2}\right)^{-\sigma}
\sigma(u_1\cosh s-u_2\sinh s)^{-1}$$
and we get 
\begin{eqnarray*}
Sf_j(a_s)
=\int_K\left(\frac{ 1-|u_1u_2^{-1}\tanh s|^2}{|1-u_1u_2^{-1}\tanh s|^2}\right)^{(\nu-\sigma)/2}
\langle \tau_j\circ \pi(k^{-1})z^j,w^j \rangle_j\\
\times \sigma(u_1\cosh s-u_2\sinh s)^{-1}\sigma(u_1)z^{\sigma}dk\\
=\int_K\left(\frac{ 1-|\tanh s|^2}{|1-u_1u_2^{-1}\tanh s|^2}\right)^{(\nu-\sigma)/2}
\langle \tau_j\circ \pi(k^{-1})z^j,w^j \rangle_j\\
\times \sigma(\cosh s-u_1^{-1}u_2\sinh s)^{-1}z^{\sigma}dk\\
=(\cosh s)^{-\nu}\int_K|1-u_1u_2^{-1}\tanh s|^{\sigma-\nu}
\langle \tau_j\circ \pi(k^{-1})z^j,w^j \rangle_j\\
\times \sigma(1-u_1^{-1}u_2\tanh s)^{-1}z^{\sigma}dk.
\end{eqnarray*}
Using the identities 
\begin{multline*}
(1-u_1^{-1}u_2\tanh s)^{-1}=\frac{1-u_2^{-1}u_1\tanh s}{|1-u_2^{-1}u_1\tanh s|^2}\\
\mbox{and}\\
|1-u_1u_2^{-1}\tanh s|=|u_1^{-1}(1-u_1u_2^{-1}\tanh s)u_1|=|1-u_2^{-1}u_1\tanh s|
\end{multline*}
in the above equality
yields
\begin{eqnarray*}
&&Sf_j(a_s)=(\cosh s)^{-\nu}\\
&&\times \int_K|1-u_2^{-1}u_1\tanh s|^{-(\sigma+\nu)}
\langle \tau_j\circ \pi(k^{-1})z^j,w^j \rangle_j\,\\
&&\qquad \qquad \qquad \qquad \times\sigma(1-u_2^{-1}u_1\tanh s)z^{\sigma}dk.
\end{eqnarray*}
We observe that the integrand is right $M$-invariant. In fact,
\begin{eqnarray*}
(\cosh s)^{-\nu}\int_K|1-u_1u_2^{-1}\tanh s|^{-(\sigma+\nu)}
\langle \tau_j\circ \pi(k^{-1})z^j,w^j \rangle_j\\
\times \sigma(1-u_2^{-1}u_1\tanh s)z^{\sigma}dk\\
=(\cosh s)^{-\nu}\int_K|1-\pi(k^{-1})^{-1}\tanh s|^{-(\sigma+\nu)}
\langle \tau_j\circ \pi(k^{-1})z^j,w^j \rangle_j\\
\times \sigma(1-\pi(k^{-1})^{-1}\tanh s)z^{\sigma}dk.
\end{eqnarray*}
Therefore, it can be written as an integral over the coset space $K/M \simeq SU(2)$, i.e., 
\begin{eqnarray*}
&&Sf_j(a_s)=(\cosh s)^{-\nu}\\
&&\times \int_{SU(2)}|1-l^{-1}\tanh s|^{-(\sigma+\nu)}
\langle \tau_j(l)z^j,w^j \rangle_j
\,\sigma(1-l^{-1}\tanh s)z^{\sigma}dl.
\end{eqnarray*}
Making the change of variables $l \mapsto l^{-1}$, and using the invariance of
the Haar measure on $SU(2)$ under this map, yields
\begin{eqnarray*}
&&Sf_j(a_s)=(\cosh s)^{-\nu}\\
&&\times \int_{SU(2)}|1-l\tanh s|^{-(\sigma+\nu)}
\langle \tau_j(l^{-1})z^j,w
^j \rangle_j
\,\sigma(1-l\tanh s)z^{\sigma}dl,
\end{eqnarray*}
and this finishes the proof.
\end{proof}

\subsection{Highest weight-vectors for K-types}
The polynomial functions $p_{l_1,l_2}$ defined by 
\begin{eqnarray*}
p_{l_1,l_2}(z,w):=z^{l_1}w^{l_2},
\end{eqnarray*}
for which $l_1+l_2=N$ 
form a basis for $V_{N}$. Occasionally we will however use the somewhat ambiguous notation $z^{l_1}w^{l_2}$
when there is no risk for misinterpretation

We write $\zeta=\tanh s$ and consider the action of $(1-\zeta l)$ on
the basis vector $p_{l_1,l_2}$. 
We have
\begin{eqnarray}
&&(1-\zeta l)p_{l_1,l_2}(z,w) \label{tensorverkan} \\
&&=\left((1-\zeta l)^{-1} \left(
\begin{array}{c}
z \\
w \\
\end{array}\right)\right)^{l_1}_1
\left((1-\zeta l)^{-1} \left(
\begin{array}{c}
z \\
w \\
\end{array}\right)\right)^{l_2}_2, \nonumber
\end{eqnarray}
where the subscripts denote the projection functions
\begin{eqnarray*}
(z,w)_1=z, \qquad (z,w)_2=w
\end{eqnarray*}
onto the first and second coordinate respectively.

The Binomial theorem gives the following expression for the first factor above:
\begin{eqnarray*}
&&\left((1-\zeta l)^{-1} \left(
\begin{array}{c}
z \\
w \\
\end{array}\right)\right)^{l_1}_1\\
&&=\det(1-\zeta l)^{-l_1}\sum_{j_1=0}^{l_1}
\left( \begin{array}{c}
l_1\\
j_1\\
\end{array}\right)
z^{j_1}(-\zeta)^{l_1-j_1}\left(l^{-1} \left(
\begin{array}{c}
z\\
w\\
\end{array}\right)\right)^{l_1-j_1}_1\\
&&=
\det(1-\zeta l)^{-l_1}\sum_{j_1=0}^{l_1}z^{j_1}
\left( \begin{array}{c}
l_1\\
j_1\\
\end{array}\right)
(-\zeta)^{l_1-j_1}
l p_{l_1-j_1,0}(z,w),
\end{eqnarray*}
and the second factor has a similar expression. Substituting these into \eqref{tensorverkan} yields the double sum
\begin{eqnarray}
&&(1-\zeta l)p_{l_1,l_2}(z,w)=\det(1-\zeta l)^{-\sigma} \\ \label{binomial}
&&\times \sum_{j_1=0}^{l_1}\sum_{j_2=0}^{l_2}\left( \begin{array}{c}
l_1\\
j_1\\
\end{array}\right)
\left( \begin{array}{c}
l_2\\
j_2\\
\end{array}\right)
(-\zeta)^{(l_1+l_2-j_1-j_2)}z^{j_1}w^{j_2}
lp_{l_1-j_1,l_2-j_2}(z,w). \nonumber
\end{eqnarray}
Denote the normalisation of the basis vector $p_{r,s}$ by $e_{r,s}$. Then $$p_{r,s}=(r!s!)^{1/2}e_{r,s}$$ and the term 
$lp_{l_1-j_1,l_2-j_2}$ in \eqref{binomial} can be written as the sum
\begin{eqnarray}
&&lp_{l_1-j_1,l_2-j_2}=(r!s!)^{1/2}\\ \label{matris}
&&\times \sum_{r+s=l_1+l_2-j_1-j_2}M(l;\,l_1-j_1, l_2-j_2;\,r,s)e_{r,s}. \nonumber
\end{eqnarray}

In what follows, we will use an expression for the first factor in the integrand in Lemma \ref{Szegö} as 
a series of $SU(2)$-characters. The following result can be found in \cite{hua}.

\begin{lemma}
The function $l \mapsto(\det (1-l \tanh s)^{-\lambda}$ has the character expansion
\begin{eqnarray*}
&&(\det (1-l \tanh s)^{-\lambda}\\
&&=\sum_{j=0}^{\infty}\sum_{i=0}^{\infty}
\frac{(\lambda-1)_{i+j+1}}{(i+j+1)!}\frac{(\lambda-1)_i}{i!}(j+1)(\tanh s)^{2i+j}\chi_j\\
&&=\sum_{j=0}^{\infty}\frac{(\lambda-1)_{j+1}}{j!}(\tanh s)^j\,\,_2F_1(\lambda+j,\lambda-1;j+2;\tanh^2s)\chi_j.
\label{characterexp}
\end{eqnarray*}
\end{lemma}

\begin{prop} \label{highestprop}
The function $Sf_j$ has the following expression when restricted to $A$.
\begin{eqnarray*}
Sf_j(a_s)=(1-\tanh^2 s)^{\frac{\sigma+2}{2}}\sum_{i=0}^{\sigma}\left( \begin{array}{c} \sigma \\i\\ \end{array}\right)
(-1)^{j+i} \frac{(2\sigma)_{j+i+1}\,j!}{(j+i+1)!}(\tanh s)^{j+2i}\\
\times_2F_1(2\sigma+1+j+i,2\sigma;j+i+2;\tanh^2s)z^{\sigma}.
\end{eqnarray*}
\end{prop}

\begin{proof}
In the defining integral $$Sf_j(a_s)=(\cosh s)^{-\nu}\int_{SU(2)}(\det(1-\zeta l))^{-(\sigma+1)}\langle \tau_j(l^{-1})z^j,w^j \rangle_j
\sigma(1-\zeta l)z^{\sigma}dl$$
we already have a character expansion for the factor $(\det(1-\zeta l))^{-\lambda}$. Therefore, it suffices to
determine the expansion of 
\begin{eqnarray*}
\langle \tau_j(l^{-1})z^j,w^j \rangle_j
\sigma(1-\zeta l)z^{\sigma}
\end{eqnarray*}
into matrix coefficients. As a special case of \eqref{binomial} we have 
$$\sigma(1-\zeta l)p_{\sigma,0}=\det(1-\zeta l)^{-\sigma}\sum_{i=0}^{\sigma}\left( \begin{array}{c} \sigma \\i\\ \end{array}\right)
(-\zeta)^{\sigma-i}z^i lp_{\sigma-i,0},$$
and the special case of \eqref{matris} is
$$lp_{\sigma-i,0}=((\sigma-i)!)^{1/2}\sum_{r=0}^{\sigma-i}\langle \tau_{\sigma-i}(l)e_{\sigma-i,0},e_{r,\sigma-i-r} 
\rangle_{\sigma-i}e_{r,\sigma-i-r}.$$
To sum up, we have
\begin{eqnarray*}
&&\sigma(1-\zeta l)p_{\sigma,0}(z,w)=\det(1-\zeta l)^{-\sigma}\\
&& \times \sum_{i=0}^{\sigma}\left( \begin{array}{c} \sigma \\i\\ \end{array}\right)
(-\zeta)^{\sigma-i}z^i((\sigma-i)!)^{1/2}\sum_{r=0}^{\sigma-i}\langle \tau_{\sigma-i}(l)e_{\sigma-i,0},e_{r,\sigma-i-r} 
\rangle_{\sigma-i}\\
&&\qquad \qquad \qquad \qquad \qquad \qquad \qquad \qquad \qquad \qquad \times e_{r,\sigma-i-r}(z,w).
\end{eqnarray*}
The integrand is thus a linear combination of terms of the form 
\begin{eqnarray}
\langle \tau_j(l^{-1})z^j,w^j \rangle_j
\langle \tau_{\sigma-i}(l)e_{\sigma-i,0},e_{r, \sigma-i-r} \rangle_{\sigma-i}. \label{matriskoeff}
\end{eqnarray}
It is easy to see (using the formula \eqref{polynomrum}) that the identity
\begin{eqnarray}
\langle \tau_j(l^{-1})z^j,w^j \rangle_j=(-1)^j\langle \tau_j(l)z^j,w^j \rangle_j
\end{eqnarray} 
holds, and using this identity in \eqref{matriskoeff}, the resulting terms are matrix coefficients for the tensor product 
representation $\tau_j \otimes \tau_{\sigma-i}$.
In fact, 
\begin{multline}
\langle \tau_j(l)z^j,w^j \rangle_j
\langle \tau_{\sigma-i}(l)e_{\sigma-i,0},e_{r, \sigma-i-r} \rangle_{\sigma-i}\\
=\left(\frac{1}{(\sigma-i)!r!(\sigma-i-r)!}\right)^{1/2}
\langle (\tau_j \otimes \tau_{\sigma-i})(l)(z^j \otimes z^{\sigma-i}), w^j \otimes z^rw^{\sigma-i-r}\rangle_{j \,\,\otimes 
\,(\sigma-i)}. \label{scalartensor}
\end{multline}
We recall the Clebsch-Gordan decomposition for the tensor product $V_j \otimes V_{\sigma-i}$. There is an isometric 
$\mbox{diag}(SU(2) \times SU(2))$-intertwining operator 
$$\phi: V_j \otimes V_{\sigma-i} \rightarrow V_{j+\sigma-i} \oplus \cdots$$ 
which is of the form $$\phi=\phi_{j+\sigma-i} \oplus \cdots \oplus \phi_{(\pm(j-(\sigma-i))}$$ where each term is 
an intertwining partial isometry and the sum is orthogonal, i.e., the terms have mutually orthogonal kernels.
The vector $z^j \otimes z^{\sigma-i}$ is a weight vector of weight $-(j+\sigma-i)\sqrt{-1}H_1^*$ and hence it maps
to a highest weight vector in the summand $V_{j+\sigma-i}$. Therefore only the term corresponding to this 
summand in the orthogonal expansion of the inner product \eqref{scalartensor} is nonzero.
To be more precise, we use the isometry $\phi$ to write the matrix coefficient \eqref{scalartensor}
as the sum
\begin{multline}
\langle (\phi(\tau_j \otimes \tau_{\sigma-i})(l)(z^j \otimes z^{\sigma-i}), \phi(w^j \otimes z^rw^{\sigma-i-r})
\rangle\\
=\sum_s \langle \phi_{j+\sigma-i-2s}(\tau_j \otimes \tau_{\sigma-i})(l)(z^j \otimes z^{\sigma-i}), 
\phi_{j+\sigma-i-2s}(w^j \otimes z^rw^{\sigma-i-r})\rangle_{j+\sigma-i-2s}. \label{scalarsum}
\end{multline}
Since $$\phi_{j+\sigma-i}((\tau_j \otimes \tau_{\sigma-i})(l)(z^j \otimes z^{\sigma-i}))
=\left(\frac{j!(\sigma-i)!}{(j+\sigma-i)!}\right)^{1/2}\tau_{j+\sigma-i}(l)z^{j+\sigma-i},$$
the sum \eqref{scalarsum} is equal to its first term 
\begin{eqnarray*}
\frac{j!(\sigma-i)!}{(j+\sigma-i)!}
\langle \tau_{j+\sigma-i}(l)z^{j+\sigma-i}, z^rw^{j+\sigma-i-r} \rangle_{j+\sigma-i}.
\end{eqnarray*}
Moreover, since we are integrating against characters, only the term corresponding to $r=j+\sigma-i$ will contribute.
Hence we have the equality
\begin{multline}
Sf_j(a_s)=(1-\tanh^2 s)^{\frac{\sigma+2}{2}}\sum_{i=0}^{\sigma}\left( \begin{array}{c} \sigma \\i\\ \end{array}\right)
(-\zeta)^{\sigma-i}(-1)^j\frac{j!(\sigma-i)!}{(j+\sigma-i)!}z^{\sigma}\\
\times \int_{SU(2)}(\det(1-\zeta l))^{-(2\sigma+1)}
\langle \tau_{j+\sigma-i}(l)z^{j+\sigma-i}, z^{j+\sigma-i} \rangle_{j+\sigma-i}dl. \label{szegösum}
\end{multline}
So, by using the character expansion \eqref{characterexp} and the Schur orthogonality relations for matrix coefficients,
we get the following expression for the above integral with the index $i$ fixed.
\begin{eqnarray*}
&&\int_{SU(2)}(\det(1-\zeta l))^{-(2\sigma+1)}
\langle \tau_{j+\sigma-i}(l)z^{j+\sigma-i}, z^{j+\sigma-i} \rangle_{j+\sigma-i}dl\\
&=&\frac{(2\sigma)_{j+\sigma-i+1}}{(j+\sigma-i)!}(\tanh s)^{j+\sigma-i}\\
&&\times 2F_1(2\sigma+1+j+\sigma-i,2\sigma;
j+\sigma-i+2;\tanh^2s)\\
&&\times \frac{(j+\sigma-i)!}{j+\sigma-i+1}.
\end{eqnarray*}
So, substitution of this into the sum \eqref{szegösum} and reversing the order of summation yields
\begin{multline*}
Sf_j(a_s)=(1-\tanh^2 s)^{\frac{\sigma+2}{2}}\sum_{i=0}^{\sigma}\left( \begin{array}{c} \sigma \\i\\ \end{array}\right)
(-1)^{j+i} \frac{(2\sigma)_{j+i+1}\,j!}{(j+i+1)!}(\tanh s)^{j+2i}\\
\times_2F_1(2\sigma+1+j+i,2\sigma;j+i+2;\tanh^2s)z^{\sigma}.
\end{multline*}
\end{proof}

We now return to the language of section 3, so that $\sigma$ corresponds to the natural number $k$.
We can now state the main theorem on the $K$-types.
\begin{thm}
For $k \geq 1$, the highest weight vector for the $K$-type $V_j^* \otimes V_{k + j}$ is the function $F_j: \mathcal{D} \rightarrow V_{\tau}$, 
whose restriction to
$A \cdot 0$ is given by
\begin{eqnarray*}
F_j(t)=(1-t^2)\sum_{i=0}^{k}\left( \begin{array}{c} \sigma \\i\\ \end{array}\right)
(-1)^{j+i} \frac{(2k)_{j+i+1}\,j!}{(j+i+1)!}t^{j+2i}\\
\times_2F_1(2k+1+j+i,2k;j+i+2;t^2)z^{k}.
\end{eqnarray*}
\end{thm}

\begin{proof}
Letting $t=\tanh s=\left( \begin{array}{cc}
\cosh s & \sinh s \\
\sinh s & \cosh s \\
\end{array}\right) \cdot 0$, and applying the trivialisation mapping \eqref{globaltrivial}, together with \eqref{realautomorf},
to the functions $Sf_j$ in Proposition \ref{highestprop} immediately
gives the result.
\end{proof}

\end{document}